\documentclass[acmsmall,screen]{MapleTrans}

\usepackage[normalem]{ulem}

\usepackage{enumitem}
\usepackage{etoolbox}
\usepackage{breqn}
\usepackage{lineno}
\usepackage{maple}
\usepackage{xfrac}
\usepackage{xcolor}
\definecolor{mygreen}{rgb}{0,0.6,0}
\definecolor{mygray}{rgb}{0.5,0.5,0.5}
\definecolor{mymauve}{rgb}{0.58,0,0.82}
\definecolor{altblue}{rgb}{0.0,0.6,1.0}
\definecolor{lstbg}{cmyk}{0.05, 0.01, 0, 0}
\definecolor{morebluish}{cmyk}{0.06,0.04,0,0}
\usepackage{listings}
\input{listings-maple-definition.sty}
\setcopyright{acmlicensed}
\copyrightyear{2024}
\acmYear{2024}
\acmDOI{10.5206/mt.vViI.XXXXX}

\acmJournal{MAPLETRANS}
\acmVolume{X}
\acmNumber{Y}
\acmArticle{Z}
\acmMonth{0}

\begin{document}

\title[MT Template 1.0]{Simple Continued Fractions an Approach for High School Students}

\author{Athanasios Paraskevopoulos}
\authornotemark[1]
\email{std156374@ac.eap.gr}
\affiliation{%
  \institution{School of Science and Technology, Hellenic Open Univeristy}
  \city{Patras}
  \country{Greece}
  \postcode{263 35}
}

\renewcommand{\shortauthors}{Paraskevopoulos}

\begin{abstract}

This paper introduces high school students to the mathematical concept of continued fractions, encompassing both finite and infinite forms. It delves into fundamental properties, the computation of quadratic numbers, and the concept of conjugate quadratic numbers. A substantial focus is placed on approximating real numbers and understanding convergence properties. By fostering an engaging and interactive learning environment, the paper aims to enhance students' mathematical proficiency and problem-solving skills. Through exploring the intricate relationships within number systems, students will gain a comprehensive understanding of continued fractions, providing a solid foundation for advanced mathematical studies.
\end{abstract}


\maketitle

\section{Introduction}
Continued fractions, a mathematical concept, are expressions derived through an iterative process of representing a number as the sum of its integer part and the reciprocal of another number. This process is repeated, with each subsequent number being expressed as the sum of its integer part and another reciprocal. In a finite continued fraction, the iteration terminates after a finite number of steps by using an integer instead of another continued fraction. Conversely, an infinite continued fraction is an unending expression. In both cases, all integers in the sequence, barring the first, must be positive. These integers, $a_i$, are referred to as the coefficients or terms of the continued fraction.

The numerator of all the fractions is typically assumed to be 1. If arbitrary values or functions are used in place of one or more numerators or integers in the denominators, the resulting expression is a generalized continued fraction. When it is necessary to distinguish the first form of generalized fractions, the former may be called simple or regular continued fractions or said to be in a canonical form.

Continued fractions exhibit numerous remarkable properties related to the Euclidean algorithm for integers or real numbers. Every rational number $\frac{p}{q}$ has two closely related expressions as a finite continued fraction, whose coefficients $a_i$ can be determined by applying the Euclidean algorithm to $(p,q)$. The numerical value of an infinite continued fraction is irrational; it is defined from its infinite sequence of integers as the limit of a sequence of values for finite continued fractions. Each finite continued fraction of the sequence is obtained by using a finite prefix of the infinite continued fraction’s defining sequence of integers. Moreover, every irrational number $\alpha$ is the value of a unique infinite regular continued fraction, whose coefficients can be found using the non-terminating version of the Euclidean algorithm applied to the incommensurable values $\alpha$ and 1. This way of expressing real numbers (rational and irrational) is called their continued fraction representation \cite{Pettofrezzo1970}.

All graphical representations and plots in this study were generated using Maple.

\section{Generalized continued fraction}

The term "continued fraction" refers to a class of expressions where a generalized continued fraction of the form
$$a_0+\frac{b_1}{a_1+\frac{b_2}{a_2+\frac{b_3}{a_3+\cdots}}}$$

(and the terms may be integers, reals, complexes, or functions of these) represents the most general variety \cite{Rockett1992}. This generalization allows for greater flexibility in mathematical modeling and analysis, as it can encompass a wider range of numerical and functional expressions. Generalized continued fractions are particularly useful in number theory and approximation theory, providing a robust framework for the representation and approximation of real numbers, complex numbers, and even functions.

The term "continued fraction" was first used by John Wallis in his \textit{Arithmetica Infinitorum} of 1653 \cite{Havil2003}, although other sources list the publication date as 1655 or 1656. An archaic term for a continued fraction is an "anthyphairetic ratio," which reflects the ancient method of reciprocation used to derive such expressions. This historical context underscores the long-standing significance of continued fractions in mathematical literature and their evolution over centuries.

Generalized continued fractions also have applications in solving equations and optimization problems. They offer a systematic approach to generating convergent sequences that approximate solutions to these problems. Additionally, they are instrumental in the study of Diophantine equations and transcendental number theory, where their unique properties facilitate the investigation of irrational numbers and their approximations.

\section{Simple Continued Fractions}

A simple continued fraction is a special case of a generalized continued fraction where the partial numerators are equal to unity, i.e., $b_n=1$ for all $n=1, 2, ....$. A simple continued fraction is thus an expression of the form

$$a_0+\frac{1}{a_1+\frac{1}{a_2+\frac{1}{a3+\cdots}}}$$

When used without qualification, the term "continued fraction" often means "simple continued fraction" or, more specifically, regular (i.e., a simple continued fraction whose partial denominators $a_0, a_1, ...$ are positive integers;\cite{Rockett1992}). The number of terms can either be finite or infinite. A more convenient way to denote continued fractions such as the one above would be to denote it by $N=[a_0,a_1,a_2,a_3,a_4,a_5,a_6,\cdots]$.

The $n$th convergent of a simple continued fraction $N=[a_0,a_1,a_2,a_3,a_4,a_5,a_6,\cdots]$ is given by the ratio of two successive terms in the sequence of numerators $P_n$ and denominators $Q_n$ defined by the recurrence relations:

\begin{align*}
P_{-1} &= 1, & Q_{-1} &= 0, \\
P_{0} &= a_0, & Q_{0} &= 1, \\
P_{n} &= a_nP_{n-1} + P_{n-2}, & Q_{n} &= a_nQ_{n-1} + Q_{n-2} \quad \text{for } n \geq 1.
\end{align*}

The sequence of convergents $\frac{P_n}{Q_n}$ provides the best rational approximations to the real number represented by the continued fraction. This property is known as the \textit{best approximation property} of continued fractions.

Every real number has a unique representation as a simple continued fraction, except for the ambiguity at the end of finite continued fractions, where $[a_0,a_1,a_2,\ldots,a_n-1,1]$ is equivalent to $[a_0,a_1,a_2,\ldots,a_n]$. The continued fraction representation of a rational number is always finite, while the continued fraction representation of an irrational number is always infinite. This property provides a method to distinguish between rational and irrational numbers.

\subsection*{Finite Continued Fractions}

A finite continued fraction is an expression such as the one shown above that could end. Every rational number can be equated to a finite continued fraction. The only skill needed would be the division of fractions.

For example, $\frac{47}{17}$ can be expressed as:

$$\frac{47}{17}=2+\frac{13}{17}=2+\frac{1}{\frac{17}{13}}=2+\frac{1}{1+\frac{4}{13}}=2+\frac{1}{1+\frac{1}{\frac{13}{4}}}=2+\frac{1}{1+\frac{1}{3+\frac{1}{4}}}$$

In the short form, $\frac{47}{17}=[2;1,3,4]$

The process of converting a rational number to a finite continued fraction involves repeated application of the Euclidean algorithm for the greatest common divisors. This algorithm is based on the observation that if $a > b$, then $\gcd(a, b) = \gcd(b, a \mod b)$.

\subsection*{Infinite Continued Fractions}

Unlike finite continued fractions, the chain of fractions never ends in an infinite continued fraction. Every irrational number can be equated to an infinite continued fraction. This fact was discovered and proven by the Swiss Mathematician, Leonhard Euler (1707-1783). Some of Euler’s infinite continued fractions are as follows:

$$\sqrt{2}=1+\frac{1}{2+\frac{1}{2+\frac{1}{2+\cdots}}}$$

A way to summarise this expression is to let $x$ be the value of the continued fraction.

$$x=1+\frac{1}{2+\frac{1}{2+\frac{1}{2+\cdots}}}\Leftrightarrow x=1+\frac{1}{1+x}$$

This equation can be solved to find $x = \sqrt{2} - 1$, which verifies Euler's representation of $\sqrt{2}$ as an infinite continued fraction.

Infinite continued fractions provide an efficient way to approximate irrational numbers. The sequence of convergents of an infinite continued fraction converges to the value of the continued fraction. The rate of convergence depends on the size of the partial quotients. For example, the continued fraction representation of $\pi$ has relatively large partial quotients, which means that the convergents of this continued fraction provide very good approximations to $\pi$.

\subsection*{Solving Quadratic Equations with Continued Fractions}

In mathematics, a quadratic equation is a polynomial equation of the second degree. The general form is $ax^{2}+bx+c=0$, where $a \neq 0$.

The quadratic equation on a number $x$ can be solved using the well-known quadratic formula, which can be derived by completing the square. That formula always gives the roots of the quadratic equation, but the solutions are expressed in a form that often involves a quadratic irrational number, which is an algebraic fraction that can be evaluated as a decimal fraction only by applying an additional root extraction algorithm.

If the roots are real, there is an alternative technique that obtains a rational approximation to one of the roots by manipulating the equation directly. The method works in many cases, and long ago it stimulated further development of the analytical theory of continued fractions \cite{Wall1948}.

Joseph-Louis Lagrange (1736-1813) proved that the continued fraction expansion of a real number $x$ is ultimately periodic, i.e,

$$x=[a_0,...a_k,b_1,...,b_h,b_1,...,b_h]$$

if and only if $x$ is a quadratic number, that is, $x$ is the root of a quadratic polynomial with rational coefficients. In such cases, we use the shorter notation

$$x=[a_0,...a_k,\Bar{b_1,...,b_h]},$$

in a way similar to how it is done for repeating decimals.

This result is of fundamental importance in the theory of continued fractions. It provides a method for finding the continued fraction expansion of a quadratic irrational number, and it also provides a method for finding the quadratic irrational number that corresponds to a given periodic continued fraction.

The process of finding the continued fraction expansion of a quadratic irrational number involves finding the period of the continued fraction, which can be done using the method of Lagrange. The process of finding the quadratic irrational number that corresponds to a given periodic continued fraction involves solving a quadratic equation, which can be done using the quadratic formula.

\subsection{Comparison with Built-in Maple Functions}
While Maple provides built-in functions for handling continued fractions, such as NumberTheory[ContinuedFraction], this paper employs custom implementations to better suit educational purposes. The custom code allows for step-by-step demonstrations, making it easier for students to understand how each convergent is calculated. Additionally, it offers more flexibility in modifying the process to explore generalized or specific continued fractions, such as those used in solving quadratic equations or demonstrating the convergence properties of irrational numbers.

\subsubsection*{Example}

Here is a simple example to illustrate the solution of a quadratic equation using continued fractions. We begin with the equation $x^2=2$ and manipulate it directly. Subtracting one from both sides we obtain $x^2-1=1$.
This is easily factored into $(x+1)(x-1)=1$ from which we obtain $ (x-1)={\frac {1}{1+x}}$ and finally $ x=1+{\frac {1}{1+x}}$.
Now comes the crucial step. We substitute this expression for $x$ back into itself, recursively, to obtain
$$ x=1+{\cfrac {1}{1+\left(1+{\cfrac {1}{1+x}}\right)}}=1+{\cfrac {1}{2+{\cfrac {1}{1+x}}}}.$$
But now we can make the same recursive substitution again, and again, and again, pushing the unknown quantity $x$ as far down and to the right as we please, and obtaining in the limit the infinite continued fraction

$$x=1+{\cfrac {1}{2+{\cfrac {1}{2+{\cfrac {1}{2+{\cfrac {1}{2+{\cfrac {1}{2+\ddots }}}}}}}}}}={\sqrt {2}}.$$
Hence the continued fraction expansion of $\sqrt2$ is given by
$$\sqrt{2}=[1,2,2,2,...]=[1,\Bar{2]}.$$

\begin{lstlisting}[]

sqrt2_actual := evalf(sqrt(2));

cf_sqrt2 := proc(n) 
 local i, cf;
  cf := 0; 
  for i from n by -1 to 1 do
   cf := 1/(cf + 2); 
 end do; 
 cf + 1;
end proc;


N := 20;

approximations := [seq([n, cf_sqrt2(n)], n = 1 .. N)];

errors := [seq([n, abs(cf_sqrt2(n) - sqrt2_actual)], n = 1 .. N)];

with(plots);

p1 := listplot(approximations, style = point, color = blue, 
symbol = diamond, labels = ["Number of Terms", "Approximation"],
view = [1 .. N, sqrt2_actual - 0.1 .. sqrt2_actual + 0.1]);

p2 := plot([sqrt2_actual, sqrt2_actual], 1 .. N, 
linestyle = dot, color = red);

display([p1, p2],
title = "Comparison of Approximations vs. Actual sqrt(2)",
legend = ["Approximations", "Actual sqrt(2)"]);

error_plot := listplot(errors, style = pointline, color = green, 
symbol = circle, labels = ["Number of Terms", "Error"], 
axis[2] = [mode = "log"], 
tickmarks = [0.1, 0.01, 0.001, 0.0001, 0.00001, 0.1e-5, 0.1e-6]],
view = [1 .. N, 0.1e-8 .. 0.1]);

display(error_plot, title = "Error in Approximation (Log Scale)");
\end{lstlisting}
\begin{figure}[ht]
    \centering
    \includegraphics[width=6cm]{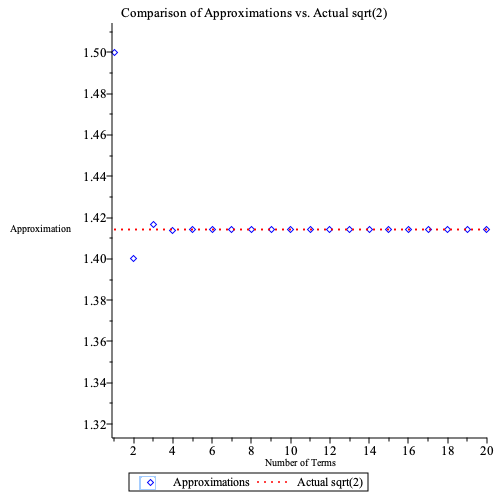}

  \caption{Convergents of Continued Fraction for $\sqrt{2} $}
\end{figure}

The first graph will provide a broader range of approximations, showing how they converge more distinctly as the number of terms increases.

\begin{figure}[ht]
    \centering
    \includegraphics[width=6cm]{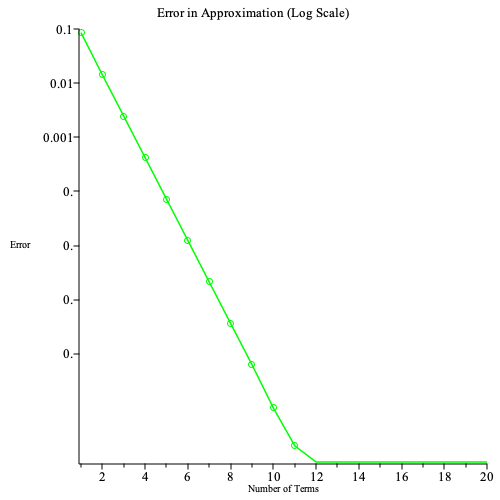}

  \caption{The Logarithmic Plot of Approximation Errors for Continued Fractions of \(\sqrt{2}\)}
\end{figure}

 This plot illustrates the convergence of continued fraction approximations to the actual value of \(\sqrt{2}\). The y-axis is displayed on a logarithmic scale to show the decreasing error across successive approximations, highlighting how the continued fractions quickly converge to the true value. Smaller error values (below \(10^{-6}\)) indicate high precision in the approximation.

By applying the fundamental recurrence formulas, we may easily compute the successive convergent of this continued fraction to be 1, 3/2, 7/5, 17/12, 41/29, 99/70, 239/169, ..., where each successive convergent is formed by taking the numerator plus the denominator of the preceding term as the denominator in the next term, then adding in the preceding denominator to form the new numerator. This sequence of denominators is a particular Lucas sequence known as the Pell numbers.

This example illustrates the power of continued fractions in solving quadratic equations and approximating irrational numbers. It also highlights the connection between continued fractions and other areas of mathematics, such as number theory and sequences.

\subsection*{Pell's Equation}

Continued fractions play an essential role in the solution of Pell's equation. For example, for positive integers $p$ and $q$, and non-square $n$, it is true that if $p^2-nq^2=\pm1$, then $\frac{p}{q}$ is a convergent of the regular continued fraction for $\sqrt{n}$. The converse holds if the period of the regular continued fraction for $\sqrt{n}$ is 1, and in general, the period describes which convergent gives solutions to Pell's equation \cite{Niven1991}.

Pell's equation is a special type of Diophantine equation, named after the English mathematician John Pell. Despite its name, the equation was first studied in detail by the ancient Indian mathematician Brahmagupta, as part of his pioneering work on number theory. The connection between Pell's equation and continued fractions was discovered much later, and it provides a powerful method for finding integer solutions to the equation.

The solutions to Pell's equation have many interesting properties and applications. For example, they can be used to approximate square roots, solve certain types of quadratic equations, and construct Pythagorean triples. The study of Pell's equation and its solutions is a fascinating topic in number theory, and it illustrates the deep connections between different areas of mathematics.

\subsection*{Task 1: Analysis of an Infinite Continued Fraction}

Consider the infinite continued fraction:

\begin{equation}
x = 1+\frac{1}{1+\frac{1}{1+\frac{1}{1+\frac{1}{1+\frac{1}{...}}}}}.
\end{equation}

This non-terminating fraction can be represented as a sequence of terms, $t_n$, defined recursively as:

\begin{equation}
t_{n+1}=1+\frac{1}{t_n},
\end{equation}

where $t_1=2$. The first few terms of this sequence, calculated to four decimal places, are:

\begin{align*}
t_1 &= 2.0 \\
t_2 &= 1.500000000 \\
t_3 &= 1.666666667 \\
t_4 &= 1.600000000\\
t_5 &= 1.625000000 \\
t_6 &\approx 1.615384615 \\
t_7 &\approx 1.619047619 \\
t_8 &\approx 1.617647059 \\
t_9 &\approx 1.618181818 \\
t_{10} &\approx 1.617977528 \\
\end{align*}

    As the graph below shows, from the $8^{th}$ term onwards, all terms converge to 1.618

\begin{lstlisting}[]
phi_actual := evalf((1 + sqrt(5))/2);

golden_ratio_cf := proc(n)
  local i, cf;
  cf := 0; 
  for i from n by -1 to 1 do 
   cf := 1/(cf + 1);
 end do;
 cf + 1; 
end proc;

N := 50;

approximations := [seq([n, golden_ratio_cf(n)], n = 1 .. N)];

errors := [seq([n, abs(golden_ratio_cf(n) - phi_actual)], n = 1 .. N)];

with(plots);

p1 := listplot(approximations, style = point, color = blue, 
symbol = diamond, labels = ["Number of Terms", "Approximation"]);

p2 := plot([phi_actual, phi_actual], 1 .. N, linestyle = dot, 
color = red);

display([p1, p2], 
title = "Comparison of Approximations vs. Actual Golden Ratio", 
legend = ["Approximations", 
"Actual Golden Ratio"]);

error_plot := listplot(errors, style = pointline, color = green, 
symbol = circle, axis[2] = [mode = "log"],
view = [1 .. N, 0.1e-8 .. 0.1],
labels = ["Number of Terms", "Error (log scale)"]);

display(error_plot, 
title = "Error in Approximation of Golden Ratio (Log Scale)");
\end{lstlisting}
\begin{figure}[ht]
    \centering
    \includegraphics[width=6cm]{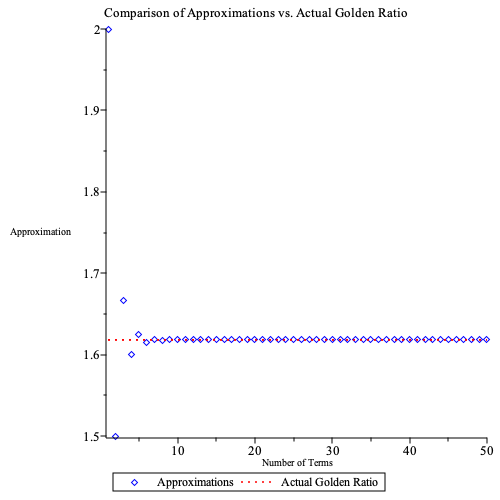}

  \caption{Sequence Converging to the Golden Ratio }
\end{figure}

In the first graph, which plots the approximations of the golden ratio versus the actual value (a red dotted line), you would expect to see that after around 10 terms, the approximations are very close to the actual value. This indicates that the continued fraction is effectively converging to the golden ratio. As the number of terms increases, the blue points (representing the approximations) should get closer and closer to the red line (the actual value), showing that the continued fraction converges well.

\begin{figure}[ht]
    \centering
    \includegraphics[width=6cm]{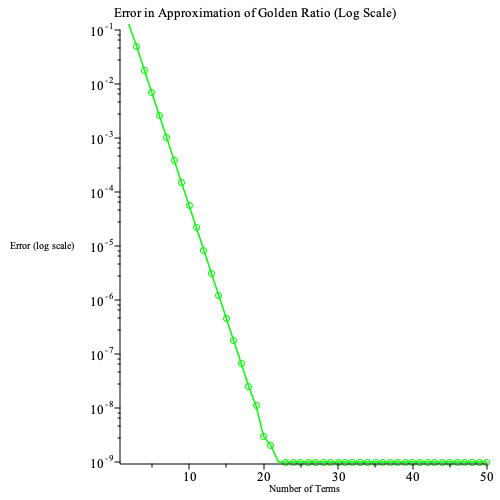}

  \caption{The Logarithmic Plot of Approximation Errors for the Golden Ratio }
\end{figure}
The second graph displays the errors of the approximations on a logarithmic scale. After around 10 terms, you should see that the errors become quite small, suggesting good convergence. By the time you reach 20 terms, the error should be very low, and the graph should demonstrate this as a steep downward trend. Since it's on a log scale, even small errors will be clearly visible, which helps highlight how effective the continued fraction is at approximating the golden ratio.
\newline

 From the $8^{th}$ term onwards, all terms converge to approximately 1.618. This limit is the golden ratio, denoted as $\varphi = \frac{1+\sqrt{5}}{2}$, which satisfies the property $\varphi = 1 + \frac{1}{\varphi}$.

By defining the sequence $R_n = \frac{t_{n+1}}{t_n}$, we can show that $\lim_{n \to \infty} R_n = \varphi$. This implies that the ratio of consecutive terms in the sequence $t_n$ converges to the golden ratio, and therefore, the infinite continued fraction can be written as:

\begin{equation}
1+\frac{1}{1+\frac{1}{1+\frac{1}{1+\frac{1}{...}}}} = \varphi.
\end{equation}

This result, which we have conjectured based on inductive reasoning, can also be proven deductively. By setting $x = 1+\frac{1}{1+\frac{1}{1+\frac{1}{1+\frac{1}{...}}}}$, we can derive the quadratic equation $x^2 - x - 1 = 0$. Solving this equation using the quadratic formula, we find that the positive solution is indeed $\varphi = \frac{1+\sqrt{5}}{2}$.

The golden ratio is a unique positive real number that appears in many different contexts, from mathematics to the arts \cite{Huntley1970}. It is the value to which the simplest infinite continued fraction converges.

\subsection*{Task 2: Analysis of a Specific Continued Fraction}

In this task, we analyze the specific continued fraction:

\begin{equation}
x = 2+\frac{1}{2+\frac{1}{2+\frac{1}{2+\frac{1}{2+\frac{1}{...}}}}}.
\end{equation}

We can represent this infinite fraction as a sequence of terms, $t_n$, where

\begin{align*}
t_1 &= 2+1, \\
t_2 &= 2+\frac{1}{3}, \\
t_3 &= 2+\frac{1}{2+\frac{1}{3}}, \\
&\vdots
\end{align*}

which can be expressed using a recursive sequence: $t_{n+1}$ is $t_{n+1}=2+\frac{1}{t_n}$.

All the terms of the sequence are finite continued fractions, so they can be computed easily. For instance,

\begin{align*}
t_1 &= 2+1 = 3.00000000, \\
t_2 &= 2+\frac{1}{3} = 2.33333333, \\
t_3 &= 2+\frac{1}{\frac{7}{3}} = 2+\frac{3}{7} = \frac{17}{7} = 2.42857143.
\end{align*}
The goal of this analysis is to generate the sequence using the recursive formula, study its convergence behavior, and compare the approximations to the known limit. The approach aligns with numerical analysis techniques that examine recursive sequences for stability and convergence. \cite{corless2013graduate}

The following table, generated by Maple, shows the first ten terms of the sequence:

\begin{table}[ht]
\centering
\begin{tabular}{|c|c|}
\hline
$n$ & $t_{n+1}=2+\frac{1}{t_n}$ \\
\hline
1 & 3.0 \\
2 & 2.333333333\\
3 & 2.428571429 \\
4 & 2.411764706 \\
5 & 2.414634146 \\
6 & 2.414141414 \\
7 & 2.414225941 \\
8 & 2.414211438 \\
9 & 2.414213927 \\
10 & 2.414213500\\
\hline
\end{tabular}
\caption{First ten terms of the sequence}
\label{tab:my_label}
\end{table}

\begin{lstlisting}[]
with(plots);
Digits := 20;

seq2 := proc(n)
 local t, i, seq; 
 t := 3.0; 
 seq := [t];
 for i to n do 
  t := 2 + 1/t; 
  seq := [op(seq), t]; 
   end do; 
   return seq; 
end proc;
seq2Values := seq2(10);

errorPlot := proc(n)
 local actualValue, t, errors, I;
 actualValue := 1 + sqrt(2); 
 t := 3.0;
 errors := [abs(actualValue - t)]; 
 for i to n do 
  t := 2 + 1/t; 
  errors := [op(errors), abs(actualValue - t)];
  end do; 
  return errors; 
end proc;

errorValues := errorPlot(10);

logErrors := [seq([i, log10(errorValues[i])], i = 1 .. 10)];

actualValue := 1 + sqrt(2);
approxPoints := [seq([i, seq2Values[i]], i = 1 .. 10)];

approxPlot := listplot(approxPoints, style = pointline, 
symbol = solidcircle, title = "Approximation Convergence to Actual Value", 
labels = ["n", "Value"], color = green);

actualLine := plot([actualValue, actualValue], 1 .. 10,
linestyle = dash, color = red);

combinedPlot := display([approxPlot, actualLine],
title = "Continued Fraction Approximation vs Actual Value", 
legend = ["Approximations", "Actual Value"], 
labels = ["Iteration (n)", "Value"]);

logPlot := listplot(logErrors, style = pointline, symbol = solidcircle,
title = "Logarithmic Error Plot for Continued Fraction Approximation",
labels = ["Iteration (n)", "Log(Error)"], color = blue);

\end{lstlisting}

\begin{figure}[ht]
    \centering
    \includegraphics[width=6cm]{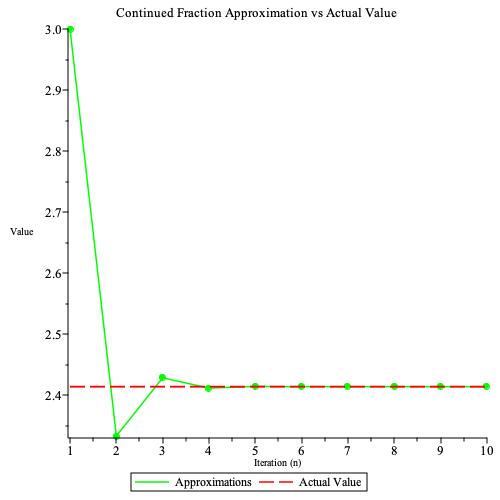}

  \caption{Sequence for $x = 2+\frac{1}{2+\frac{1}{2+\frac{1}{2+\frac{1}{2+\frac{1}{...}}}}}$}
\end{figure}
The code generates a sequence of 10 terms using the recursive formula \(t_{n+1} = 2 + \frac{1}{t_n}\) to approximate the continued fraction for \(1 + \sqrt{2}\). The procedure begins from an initial term \(t_0 = 3.0\) and produces a visual plot showing how the approximations converge to the expected value.

\begin{figure}[ht]
    \centering
    \includegraphics[width=6cm]{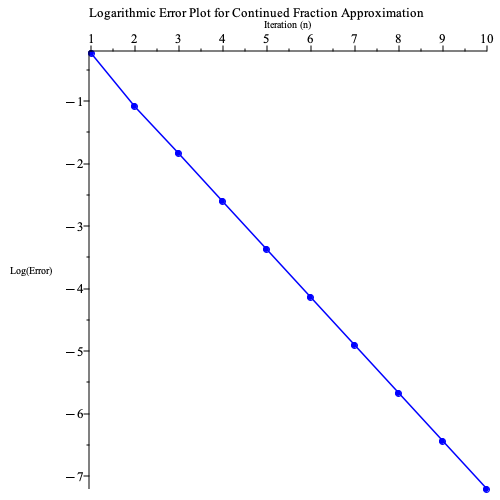}

  \caption{Logarithmic Error Plot }
\end{figure}
The `errorPlot` procedure calculates the absolute differences between each term and the known limit, and the logarithmic transformation reveals the rate of error reduction over iterations. The plot displays how quickly the approximations become accurate. The rapid decline in error suggests efficient convergence.\\
  
We observe that, as $n$ becomes very large, the terms tend to 2.414.\\

If we want to determine, for example, the 200th term, we face the problem that we have to calculate all the previous terms until we reach the 200th, which is time-consuming. From the table above, we could assume that since the terms from the 4th onwards are at 2.414, the 200th term also tends to 2.414. However, this is not proof.\\

To prove this, we can let $x$ be the infinite continued fraction:

\begin{equation}
x = 2+\frac{1}{2+\frac{1}{2+\frac{1}{2+\frac{1}{2+\frac{1}{...}}}}}.
\end{equation}

Then, it must be true that $x = 2+\frac{1}{x}$, because of the property that infinity has, we can add another level to it any time we want. It will not change it. If we multiply both sides by $x$, we get $x^2 = 2x+1$. This is a quadratic equation that can be put in standard form as $x^2 - 2x - 1 = 0$. Using the quadratic formula, we find that

\begin{equation}
x = \frac{-(-2) \pm \sqrt{(-2)^2 - 4(1)(-1)}}{2(1)} = \frac{2 + 2\sqrt{2}}{2} = 1 + \sqrt{2} \approx 2.414.
\end{equation}

Since $x$ is a string of positive terms, it is not going to be $1 - \sqrt{2}$. Therefore, the 200th term also tends to be 2.414.
  
   \newpage
\subsection*{Task 3 and 4: Analysis of General Continued Fractions}

General continued fractions can be represented using a recursive approach, where the recursive formula takes the form:

\begin{equation}
x = \kappa+\frac{1}{\kappa+\frac{1}{\kappa+\frac{1}{\kappa+\frac{1}{\kappa+\frac{1}{...}}}}}.
\end{equation}

with different values of $\kappa$ leading to various types of sequences. In this task, we study how different values of $\kappa$ affect the convergence of the sequence, with a focus on the recursive behavior, stability, and convergence patterns. \cite{corless1992continued}
Just as the two other continued fractions we have considered, this general continued fraction can be represented as a sequence of terms, $t_n$, where

\begin{align*}
t_1 &= \kappa+1, \\
t_{n+1} &= \kappa+\frac{1}{t_n}.
\end{align*}

We construct tables, like in the previous cases for $\kappa=-3$, $\kappa=\frac{2}{3}$, and $\kappa=\sqrt{5}$.

\begin{table}[ht]
\centering
\begin{tabular}{|c|c|c|c|}
\hline
$n$ & $t_{n+1}=-3+\frac{1}{t_n}$ & $t_{n+1}=\frac{2}{3}+\frac{1}{t_n}$ & $t_{n+1}=\sqrt{5}+\frac{1}{t_n}$ \\
\hline
1 & -2. & 1.6667 & 3.236 \\
2 & -3.500000000 & 1.266654667 & 2.545091463 \\
3 & -3.285714286 & 1.456147830 & 2.628981165 \\
4 & -3.304347826 & 1.353410127 & 2.616443467\\
5 & -3.302631579 & 1.405540996 & 2.618266183 \\
6 & -3.302788845 & 1.378136487 & 2.618000115 \\
7 & -3.302774427 & 1.392284201& 2.618038930 \\
8 & -3.302775749 & 1.384910829 & 2.618033267 \\
9 & -3.302775628 & 1.388734816 & 2.618034094 \\
10 & -3.302775639 & 1.386746547 & 2.618033973 \\
\hline
\end{tabular}
\caption{First ten terms of the sequence for different values of $\kappa$}
\label{tab:my_label1}
\end{table}

We implemented a Maple procedure, generalSeq, to generate the sequence by iteratively applying this recursive formula. Starting from an initial value for example $t=-2$, the sequence approximates the continued fraction by generating a series of terms that progressively approach the limit. The method draws from established recursive algorithms in symbolic computation, as discussed in works by Corless and others \cite{corless2013graduate}.\\

Below is a detailed write-up in which we can explain the use of the `generalSeq` and `generalErrorPlot` procedures, the choice of \(\kappa\) and \(t_0\), and the resulting plots. \\

The recursive formula \(t_{n+1} = \kappa + \frac{1}{t_n}\) can produce sequences that either converge to a limit or exhibit other patterns depending on the chosen \(\kappa\) and \(t_0\). To explore this, we implemented a Maple procedure `generalSeq`, which generates the sequence by iteratively applying the recursive formula for a specified number of iterations \(n\). The sequence can be visualized to understand how it evolves over time.\\

To quantify the convergence, we introduced a second procedure, `generalErrorPlot`, which calculates the absolute error between each term of the sequence and an expected limit \(L\). This allows us to observe how closely the sequence approaches the limit over successive iterations. The errors are plotted on a logarithmic scale, making it easier to detect trends in convergence.

We analyzed the recursive sequences for the following values of \(\kappa\):
\begin{enumerate}

    \item \(\kappa = -3\): Initial value \(t_0 = -2\). The expected limit was approximately \(L = -3.302775\). 
\item \(\kappa = \frac{2}{3}\): Initial value \(t_0 = 1.6667\). The expected limit was \(L = 1.5\).
\item \(\kappa = \sqrt{5}\): Initial value \(t_0 = 3.236\). The expected limit was approximately \(L = 2.618\).

\end{enumerate}

If this infinite continued fraction converges at all, it must converge to one of the roots of the monic polynomial $x^2+\kappa x+1=0$. That depends on both the coefficient $\kappa$ and the value of the discriminant, $b^2-4ac$.

\begin{lstlisting}[]
with(plots);

generalSeq := proc(kappa, t0, n)
  local seq, i, t; 
  t := evalf(t0); 
  seq := [t]; 
  for i to n do 
   t := evalf(kappa + 1/t);
   seq := [op(seq), t]; 
  end do;
  return seq;
end proc;

generalErrorPlot := proc(kappa, t0, n, expectedValue) 
  local seq, errors, i, t; 
  t := evalf(t0); 
  seq := [t];
  errors := [abs(expectedValue - t)]; 
  for i to n do
   t := evalf(kappa + 1/t); seq := [op(seq), t]; 
   errors := [op(errors), abs(expectedValue - t)]; 
  end do;
  return errors; 
end proc;

sequences := [generalSeq(-3, -2, 9), 
generalSeq(2/3, 1.6667, 9),
generalSeq(sqrt(5), 3.236, 9)];

plots:-display(
[
plots:-listplot(sequences[1], style = pointline, 
color = green, symbol = solidcircle), 
plots:-listplot(sequences[2], style = pointline, 
color = purple, symbol = solidcircle),
plots:-listplot(sequences[3], style = pointline, 
color = orange, symbol = solidcircle)], 
legend = ["&kappa; = -3", "&kappa; = 2/3", "&kappa; = &Sqrt;5"], 
labels = ["n", "Value"], gridlines = true,
title = "Sequences for Different Values of &kappa;");

errorValues_kappa_neg3 := generalErrorPlot(-3, -2, 20, -3.302775);
logErrors_kappa_neg3 := [seq([i, log10(errorValues_kappa_neg3[I])], 
i = 1 .. 20)];

logPlot_kappa_neg3 := plots:-listplot(logErrors_kappa_neg3,
style = pointline, 
symbol = solidcircle, 
title = "Logarithmic Error Plot for &kappa; = -3", 
labels = ["Iteration (n)", "Log(Error)"], color = green);

errorValues_kappa_2over3 := generalErrorPlot(2/3, 1.6667, 20, 1.5);

logErrors_kappa_2over3 := [seq([i, log10(errorValues_kappa_2over3[I])], 
i = 1 .. 20)];

logPlot_kappa_2over3 := plots:-listplot(logErrors_kappa_2over3,
style = pointline, 
symbol = solidcircle, 
title = "Logarithmic Error Plot for &kappa; = 2/3", 
labels = ["Iteration (n)", "Log(Error)"], color = purple);

errorValues_sqrt5 := generalErrorPlot(1 + sqrt(2), 3.236, 20, 2.618);

logErrors_sqrt5 := [seq([i, log10(errorValues_sqrt5[i])], i = 1 .. 20)];

logPlot_sqrt5 := plots:-listplot(logErrors_sqrt5, 
style = pointline, 
symbol = solidcircle, 
title = "Logarithmic Error Plot for &kappa; = &Sqrt;5", 
labels = ["Iteration (n)", "Log(Error)"], color = orange);



\end{lstlisting}
  \begin{figure}[ht]
    \centering
    \includegraphics[width=6cm]{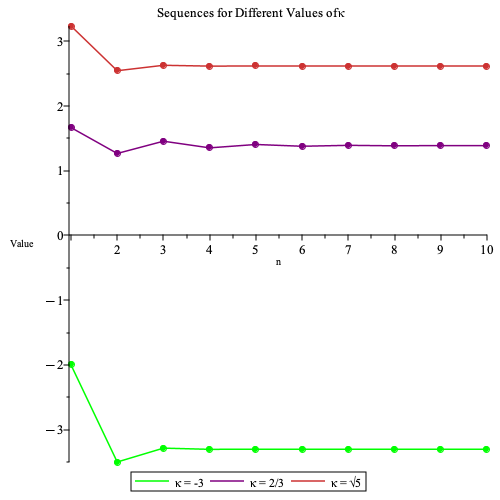}
  \caption{Sequences for Different Values of $\kappa$}
\end{figure}

This visual comparison allows us to see how quickly each sequence stabilizes, depending on the choice of \(\kappa\) and the initial term \(t_0\).\\

These values were chosen to illustrate different convergence behaviors, as recursive sequences can show varying levels of stability depending on $\kappa$.\\

The green plot represents \(\kappa = -3\), showing convergence towards a negative value.
The purple plot represents \(\kappa = \frac{2}{3}\), converging to a rational number.
The orange plot represents \(\kappa = \sqrt{5}\), converging to an irrational number.\\
\newpage
\subsubsection*{Logarithmic Error Analysis}

To further analyze the convergence, we plotted the errors on a logarithmic scale. The error \(E_n\) for each term \(t_n\) was calculated as:
\[
E_n = |t_n - L|,
\]
where \(L\) is the expected limit. 
The errors were plotted on a logarithmic scale to highlight the rate of convergence. This approach is commonly used in numerical analysis to evaluate the precision and stability of approximations.\cite{corless1992continued}

The logarithmic error plots (Fig. 8,9,10) below show how the error diminishes across 20 iterations, providing insight into the rate of convergence.\\

Green Plot (\(\kappa = -3\)): The error plot shows a steady decline, indicating convergence to the negative limit.

\begin{figure}[ht]
    \centering
    \includegraphics[width=6cm]{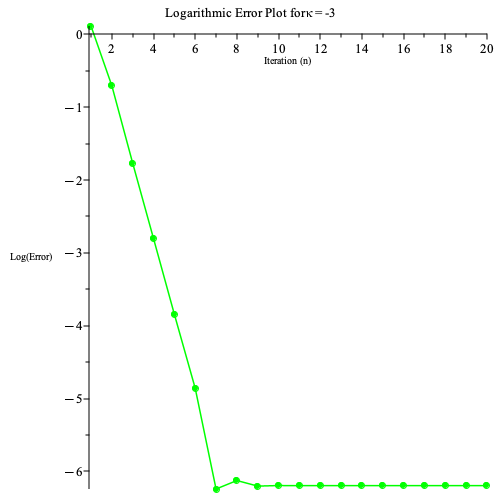}
  \caption{Logarithmic Error Plot for $\kappa=-3$}
\end{figure}
\newpage
Purple Plot (\(\kappa = \frac{2}{3}\)): Rapid reduction in error, highlighting quick convergence to the rational limit.
\begin{figure}[ht]
    \centering
    \includegraphics[width=6cm]{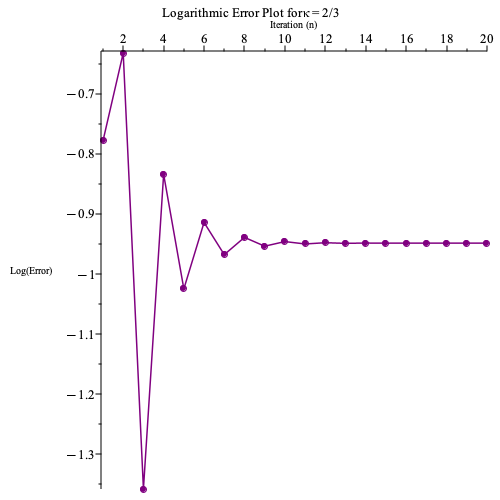}
  \caption{Logarithmic Error Plot for $\kappa=2/3$}
\end{figure}

Orange Plot (\(\kappa = \sqrt{5}\)): Gradual convergence, suggesting the sequence is approaching an irrational limit.

\begin{figure}[ht]
    \centering
    \includegraphics[width=6cm]{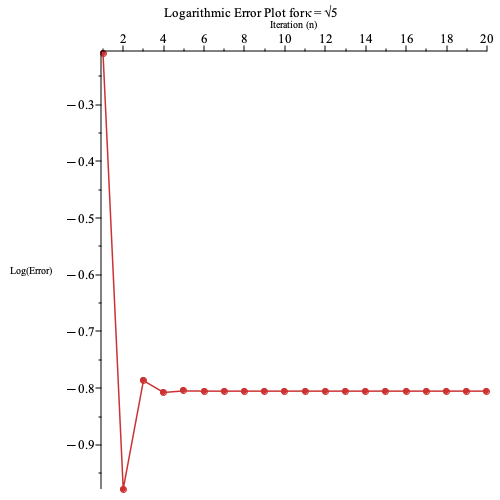}
  \caption{Logarithmic Error Plot for $\kappa=\sqrt{5}$}
\end{figure}

The logarithmic scale effectively reveals the rate at which the errors decrease, demonstrating different behaviors depending on \(\kappa\). For \(\kappa = \frac{2}{3}\), we see faster convergence compared to \(\kappa = \sqrt{5}\), where the error reduces more slowly.

This analysis shows that the behavior of the recursive sequence \(t_{n+1} = \kappa + \frac{1}{t_n}\) depends significantly on the choice of \(\kappa\). By visualizing both the sequence evolution and the error convergence, we can gain a deeper understanding of how quickly the sequence stabilizes and approaches the expected limit. 

For negative values of \(\kappa\) (e.g., \(\kappa = -3\)), the sequence can still converge, but the limit is negative. Rational and irrational values of \(\kappa\) lead to different convergence behaviors, as seen with \(\kappa = \frac{2}{3}\) and \(\kappa = \sqrt{5}\).

Assuming the existence of the limit, we have

\begin{equation}
\lim_{n\rightarrow\infty} t_n = \lim_{n\rightarrow\infty} \left(\kappa+\frac{1}{t_n}\right) = \kappa + \lim_{n\rightarrow\infty} \frac{1}{t_n}.
\end{equation}

This leads to the equation $x = \kappa + \frac{1}{x}$, which simplifies to the quadratic equation $x^2 - \kappa x - 1 = 0$.

In general, the above result could be used as in the previous cases to compute square roots to a given accuracy.

As we have already shown, continued fractions are most conveniently applied to solve the general quadratic equation expressed in the form of a monic polynomial $x^2+bx+c=0$. Starting from this monic equation we see that $$x^2- \kappa x=1 \Leftrightarrow x-\kappa=\frac{1}{x}\Leftrightarrow x=\kappa+\frac{1}{x}$$

But now we can apply the last equation to itself recursively to obtain

$$x=\kappa+\frac{1}{\kappa+\frac{1}{\kappa+\frac{1}{\kappa+\frac{1}{\kappa+\frac{1}{...}}}}}$$

If this infinite continued fraction converges at all, it must converge to one of the roots of the monic polynomial $x^2 + \kappa x + 1 = 0$. That depends on both the coefficient $\kappa$ and the value of the discriminant, $b^2-4ac$.

In general, by applying a result obtained by Euler in 1748 it can be shown that the continued fraction solution to the general monic quadratic equation with real coefficients

$$x^2+bx+c=0$$

which can always be obtained by dividing the original equation by its leading coefficient. Starting from this monic equation we see that

$$x^2+b x=-c \Leftrightarrow x+b=-\frac{c}{x}\Leftrightarrow x=-b-\frac{c}{x}$$

But now we can apply the last equation to itself recursively to obtain

$$x=-b-\frac{c}{-b-\frac{c}{-b-\frac{c}{-b-\frac{c}{-b-\frac{c}{...}}}}}$$

either converges or diverges depending on both the coefficient b and the value of the discriminant, $b^2-4ac$.

If $b = 0$ the general continued fraction solution is totally divergent; the convergent alternate between 0 and $\infty$.  If $b \neq 0$ we distinguish three cases. 

\begin{enumerate}
    \item If the discriminant is negative, the fraction diverges by oscillation, which means that its convergent wanders around in a regular or even chaotic fashion, never approaching a finite limit.
    \item If the discriminant is zero, the fraction converges to the single root of multiplicity two.
    \item If the discriminant is positive the equation has two real roots, and the continued fraction converges to the larger (in absolute value) of these. The rate of convergence depends on the absolute value of the ratio between the two roots: the farther that ratio is from unity, the more quickly the continued fraction converges.
\end{enumerate}

\section{The history of continued fractions}
 \begin{enumerate}
     \item 300 BCE Euclid’s Elements contains an algorithm for the greatest common divisor, whose modern version generates a continued fraction as the sequence of quotients of successive Euclidean divisions that occur in it.
     \item	499 The Aryabhatiya contains the solution of indeterminate equations using continued fractions.
     \item	1572 Rafael Bombelli, L’Algebra Opera – method for the extraction of square roots which is related to continued fractions.
     \item	1613 Pietro Cataldi, Trattato del modo di trovar la ratice quadra delli numeri—first notation for continued fractions.  Cataldi represented a continued fraction as $a_0 \& \frac{n_1}{d_1.} \& \frac{n_2}{d_2.} \& \frac{n_3}{d_3.}$  with dots indicating where the following fractions went.
     \item 1695 John Wallis, Opera Mathematica – introduction of the term “continued fraction”.
     \item 	1737 Leonhard Euler, De fractionibus continuis dissertation — Provided the first then-comprehensive account of the properties of continued fractions and included the first proof that the number $e$  is irrational. \cite{Sandifer2006}
     \item 1748 Euler, Introductio in analsin infinitorum. Vol.I, Chapter 18 -- proved the equivalence of a certain form of continued fraction and a generalized infinite series, proves that every rational number can be written as a finite continued fraction, and proved that the continued fraction of an irrational number is infinite. \cite{Euler1748}
     \item 	1761 Johann Lambert – gave the first proof of the irrationality of $\pi$  using a continued fraction for $\tan(x)$.
     \item 	1768 Joseph-Louis Lagrange – provide the general solution to Pell’s equation using continued fraction’s like Bombelli’s
     \item 	1770 Lagrange -- proved that quadratic irrationals expand to periodic continued fractions.
     \item 1770 Lagrange -- proved that quadratic irrationals expand to periodic continued fractions.
     \item 	1813 Carl Friedrich Gauss, Werke, Vol.3, pp. 134-138 -- derived a very general complex-valued continued fraction via a clever identity involving the hypergeometric function.
     \item 	1892 Henri Pade defined Pade approximant.
     \item 1892 Henri Pade defined Pade approximant.
     \item 1972 Bill Gosper -- First exact algorithms for continued fraction arithmetic.
     
 \end{enumerate}

\section{Conclusion}
Continued fractions serve as a versatile and powerful mathematical tool, offering unique representations of numbers and solutions to complex problems. Their utility extends particularly into number theory, where they facilitate solving equations, proving theorems, and exploring numerical properties. In this paper, we have provided an introduction to continued fractions, explored their fundamental properties, and demonstrated their applications in solving quadratic equations. Additionally, we have highlighted their connections to the Fibonacci sequence and the golden ratio. By examining the rich history of continued fractions, we underscore the significant contributions of various mathematicians over the centuries. We hope that this paper will ignite a deeper interest in continued fractions and inspire further exploration into this captivating area of mathematics.

\section{Acknowledgments}

I would like to thank my students for their curiosity and enthusiasm, which inspired me to explore new ways of teaching and learning mathematics. I would also like to thank my colleagues at the Hellenic Open University for their support and encouragement.


\bibliographystyle{plain}
\bibliography{sample-base}

\appendix

\end{document}